\documentclass[matann]{svjour}

\usepackage{amsmath}
\usepackage{amssymb}
\usepackage{amscd}

\newcommand{\gothic}{\mathfrak}
\newcommand{\p}{{\gothic{p}}}

\newcommand{\m}{{\gothic{m}}}
\newcommand{\n}{{\gothic{n}}}
\newcommand{\Spec}{\operatorname{Spec}}

\newcommand{\Hom}{\operatorname{Hom{}}}
\newcommand{\Ext}{\operatorname{Ext{}}}
\newcommand{\Tor}{\operatorname{Tor{}}}
\newcommand{\rank}{\operatorname{rank{}}}
\newcommand{\depth}{\operatorname{depth}}

\newcommand{\type}{\operatorname{type}}
\newcommand{\syz}{\operatorname{syz{}}}
\newcommand{\Ann}{\operatorname{Ann{}}}
\renewcommand{\hat}{\widehat}
\renewcommand{\phi}{\varphi}
\renewcommand{\to}{{\longrightarrow}}
\newcommand{\xx}{{\underline{x}}}

\newcommand{\limq}{{\displaystyle\lim_{q\rightarrow\infty}}}

\spnewtheorem{thm}{Theorem}[document]{\bf}{\it}
\spnewtheorem{cor}[thm]{Corollary}{\bf}{\it}
\spnewtheorem{lem}[thm]{Lemma}{\bf}{\it}
\spnewtheorem{defn}[thm]{Definition}{\bf}{\it}
\spnewtheorem{eg}[thm]{Example}{\bf}{\rm}
\spnewtheorem{prop}[thm]{Proposition}{\bf}{\it} 
\spnewtheorem*{rem}{Remark}{\bf}{\rm}

\begin{document}

\title{Two theorems about maximal Cohen--Macaulay modules}

\titlerunning{Two Theorems about MCM modules}

\author{Craig Huneke \and Graham J. Leuschke
        \thanks{Both authors were partially supported by the National
Science
        Foundation.  The second author was also partially supported by
the
        Clay Mathematics Institute.}}

\institute{Department of Mathematics \\ 
        University of Kansas \\ 
        Lawrence, KS
        66045\\ 
        \email{huneke@math.ukans.edu \and gleuschke@math.ukans.edu}} 

\date{\today}

\maketitle

\bibliographystyle{amsplain}

\begin{abstract} This paper contains two theorems concerning the
theory of maximal Cohen--Macaulay modules. The first theorem
proves that certain Ext groups between maximal Cohen--Macaulay
modules $M$ and $N$ must have finite length, provided only finitely many
isomorphism 
classes of maximal Cohen--Macaulay modules exist having ranks up to the
sum of
the ranks of $M$ and $N$. This has several corollaries. In
particular
it proves that a Cohen--Macaulay local ring of finite Cohen--Macaulay
type has an isolated
singularity. A well-known theorem of Auslander gives the same conclusion
but requires that the ring be Henselian.
Other corollaries of our result  include statements concerning when a
ring is
Gorenstein or
a complete intersection on the punctured spectrum, and the recent
theorem of Leuschke and Wiegand that the completion of an excellent
Cohen--Macaulay local ring of finite 
Cohen--Macaulay type is again of finite Cohen--Macaulay type. The second
theorem proves that a 
complete local Gorenstein domain of positive characteristic $p$ and
dimension
$d$  
is $F$-rational if and only if the number of copies of $R$ splitting out
of $R^{1/p^e}$ divided by $p^{de}$ has a positive limit. This result
generalizes
work of Smith and Van den Bergh. We call this limit the 
$F$-signature of the ring and give some of its properties.
\end{abstract}

Throughout this paper we will work with a Noetherian local ring
$(R,\m)$.
Recall that a finitely generated $R$-module $M$ is {\it maximal
Cohen--Macaulay (MCM)} if $\depth M = \dim R$.
We say that $R$ has {\it finite
Cohen--Macaulay type} (or {\it finite CM type}) provided there are only
finitely many isomorphism classes of indecomposable MCM
$R$-modules. There is a large body of work devoted to the classification
of Cohen--Macaulay local rings of finite CM type, for example see the
book \cite{Yoshino:book}. One of the first and most famous results is
that of M.~Auslander \cite{Auslander:isolsing}: if $R$ is complete
Cohen--Macaulay 
and has
finite
Cohen--Macaulay type, then $R$ has an isolated singularity, i.e. for
all primes $\p\ne \m$, $R_\p$ is a regular local ring.  (Yoshino
\cite{Yoshino:book} points
out
that $R$ need only be Henselian and have a canonical module.)   A key
point
in Auslander's proof 
is to prove that the modules Ext$^1_R(M,N)$ are of finite length for
arbitrary MCM modules $M$ and $N$, and he accomplishes this by using the
theory of 
almost split sequences.  The first result in this paper gives a
proof of the finite length of the Ext modules but with several notable
improvements: we do not need to 
assume the base ring is Henselian (nor has a canonical module), and
furthermore we can give 
an explicit bound on the power of the maximal ideal annihilating
$\Ext^1_R(M,N)$ in terms of the number of isomorphism
classes of MCM modules of multiplicity at most the sum of the
multiplicities
of $M$ and $N$ (without assuming there are only finitely many
isomorphism
classes of MCMs). This result has several corollaries of interest.
By applying the theorem to special modules $M$
we are able to give several
results which have the flavor that the ring must be either
Gorenstein, a complete intersection, or regular on the punctured
spectrum provided there are only finitely many isomorphism classes
of MCM modules up to some number depending upon the number of generators
of the special module (e.g., when the module $M$ is the canonical
module,
we  conclude the ring is Gorenstein on
the punctured spectrum).
We are able to give a direct proof
of a 
recent theorem of 
Leuschke and Wiegand that a Cohen--Macaulay local ring $R$ has finite
Cohen--Macaulay
type
if 
and only the completion of $R$ has finite Cohen--Macaulay type.

In the second theorem we discuss a more general class of rings
introduced in a paper of K.E. Smith and M. Van den Bergh
\cite{Smith-vandenBergh:1997}.
We consider reduced Cohen--Macaulay local rings of
prime
characteristic $p$, and let $q = p^e$ be a varying power of $p$.
We say that $R$ is \it $F$-finite \rm if the  Frobenius map
$F:R\rightarrow R$
sending $r$ to $r^p$ is a finite map.
A reduced local $F$-finite  Cohen--Macaulay ring $R$ is said to have
finite
$F$-representation type (or FFRT for short) if  only finitely many
isomorphism classes of indecomposable MCM modules occur as direct
summands
of $R^{1/q}$ as $q$ ranges over all powers of $p$. 
Among other results, Smith and  Van den Bergh  proved that if
we denote by $a_q$ the number of copies of $R$ splitting out of
$R^{1/q}$, then provided $R$ is strongly $F$-regular (a tight  closure
notion --- see \cite{Huneke:CBMS}), and has FFRT, the limit of
$a_q/q^{\dim R}$
exists and is positive.  We show (without the assumption of FFRT) that
if the
limit is positive, then $R$ is 
weakly $F$-regular, that is, all ideals of $R$ are tightly closed.
When $R$ is Gorenstein, we are able to prove the limit
always exists, and moreover prove that this limit
is positive if and only if $R$ is $F$-rational (which is equivalent to
strongly $F$-regular in the $F$-finite Gorenstein case). When the limit
exists, we call it the \it $F$-signature \rm of $R$ and denote
it by $s(R)$. Several examples at the end of this paper show that
$s(R)$ is a delicate invariant of the ring which gives considerable
information about the type of the singularity of $R$.  For example,
when $R$ is the hypersurface $x^2+y^3+z^5 = 0$ (the famous $E_8$
singularity), then the $F$-signature is $1/120$. (This ring
is the ring of invariants of a group of order $120$ acting on 
$k[x,y,z]$.)
We also prove that $0\leq s(R)\leq 1$, and $s(R) = 1$ if and only if the
ring is regular.

\section{ Finite type }

\begin{thm}\label{main} Let $(R,\m)$ be a local ring and let
$M$ and $N$ be two finitely generated MCM modules over $R$,
having multiplicities $m$ and $n$ respectively. Assume that there are
only finitely many isomorphism classes of MCM modules of
multiplicity
$m+n$.  If $h$ is the
number of such isomorphism classes, then $\m^h$ annihilates 
$\Ext^1_R(M,N)$. In particular, $\Ext^1_R(M,N)$ has finite length.
\end{thm}

\begin{proof} 
We claim that for any $\chi \in\Ext_R^1(M,N)$ and any $r_1,...,r_h \in
\m$,
$r_1\cdots r_h\chi =
0$.  Let
$$\CD \chi:0 @>>> N @>>> K @>>> M @>>> 0\endCD$$
be given, and consider
$$\CD r_1\cdots r_k\chi:0 @>>> N @>>> K_k @>>> M @>>> 0\endCD$$
as $k$ runs through all positive integers, and each $r_i \in \m$.  Since
each
$K_k$ is a MCM module and the multiplicity
of $K_k$ is equal to the sum of the multiplicities of $M$ and $N$, by
assumption there must be repetitions among the $K_k$.  That
is, there exist integers $a$ and $b$, with $a < b$ and $a\leq h$, such
that
$K_a \cong K_b$.
Replace $\chi$ by $r_1\cdots r_a \chi$ and set $r = r_{a+1}\cdots r_b$
to assume that we have
$$\CD \chi:0 @>>> N @>>> K @>>> M @>>> 0\endCD$$
and
$$\CD r\chi:0 @>>> N @>>> L @>>> M @>>> 0\endCD$$
and that $L\cong K$.  We will show that $\chi=0$.

Recall that $r\chi$ is constructed from $\chi$ by the following pushout
diagram
$$\CD \chi:0 @>>> N @>i>> K @>>> M @>>> 0 \\
        && @VrVV  @VVV  @| \\
        r\chi: 0 @>>> N @>>> L @>>> M @>>> 0. \endCD$$
In particular, $L \cong N\oplus K/ \langle(rn, i(n)) : n \in
N\rangle$.  This
gives another exact sequence
$$\CD \zeta:0 @>>> N
@>\left[\smallmatrix{r}\\{i}\endsmallmatrix\right]>> N\oplus K @>>> L
@>>> 0.\endCD$$
Since $L\cong K$, $\zeta$ is ``apparently split'', and
Miyata's theorem \cite[A3.29]{Eisenbud:book} implies that $\zeta$
splits.  The map ${\zeta}^*: \Ext^1_R(M,N) \to \Ext^1_R(M,N\oplus K)$
obtained by applying $\Hom_R(M,-)$ to $\zeta$ is thus a split
injection.  We claim that $\chi$ goes to zero in the second component
$\Ext^1_R(M,K)$.  Applying $\Hom_R(M,-)$ to $\chi$ gives the exact
sequence 
$$\CD \Hom_R(M,M) @>>>\Ext_R^1(M,N)  @>>>\Ext^1_R(M,K).\endCD$$
The identity endomorphism of $M$ maps to $\chi \in \Ext^1_R(M,N)$, so
$\chi$ goes to zero in $\Ext^1_R(M,K)$, and it is easy to check that
this map coincides with the second component of $\zeta^*$.  This proves
the claim, that is, $\zeta^*(\chi) = (r\chi,0)$ in $\Ext^1_R(M,N)\oplus
\Ext^1_R(M,K)$.  Letting $f$ be a left splitting for $\zeta^*$, we see
that $\chi=f(r\chi,0)=rf(\chi,0)$.  Iterating shows that $\chi$ is
infinitely divisible by $r$, so $\chi=0$, as desired.
Since $r_1,...,r_a$ are
arbitrary elements in $\m$, this proves that $\m^a$ annihilates
$\Ext^1_R(M,N)$. \hfill\qed\end{proof}

Our first application is to local rings of finite Cohen--Macaulay type.  

\begin{cor}\label{maincor} Let $(R, \m)$ be a Cohen--Macaulay local ring
of
finite
CM
type. Then $R$ is an isolated singularity.\end{cor}

\begin{proof} Set $d = \text{dim}(R)$. It suffices to prove that
$\Ext^1_R(M,N)$
has finite length for all MCM modules $M$ and $N$. For suppose this
is true, and let $\p \in \Spec R$ be different from $\m$.  Consider
syzygies
of $R/\p$ :
\begin{equation}\label{MCMsyzygies}
0 \to \syz^R_{d+1}(R/\p) \to F_{d} \to \syz^R_{d}(R/\p) \to 0.
\end{equation}
Since $R$ is Cohen--Macaulay, $M:=\syz^R_{d}(R/\p)$ and
$N:=\syz^R_{d+1}(R/\p)$ are MCM $R$-modules.  By hypothesis, localizing
(\ref{MCMsyzygies}) at $\p$ gives a split exact sequence, so $M_\p$ and
$N_\p$ are both free $R_\p$-modules.  But this
implies that the residue field of $R_\p$ has a finite free resolution,
and so $R_\p$ is a regular local ring.

 Since $R$ has finite CM type, the assumptions of Theorem~\ref{main} are
satisfied for all such $M$ and $N$, giving the desired conclusion.
\hfill\qed\end{proof}

\begin{rem}Corollary~\ref{maincor} is a generalization of a celebrated
theorem of Auslander \cite{Auslander:isolsing}, which has the same
conclusion but requires that $R$ be complete.  In fact, Auslander
obtained this theorem as a corollary of a more general result, which we
can also recover.  Specifically, Auslander considers the following
situation: Let $T$ be a complete regular local ring and let $\Lambda$ be
a possibly noncommutative $T$-algebra which is a finitely generated free
$T$-module.  Say that $\Lambda$ is {\it nonsingular} if
$\operatorname{gl.dim.}\Lambda = \dim T$, and that $\Lambda$ has {\it
finite representation type} if there are only finitely many isomorphism
classes of indecomposable finitely generated $\Lambda$-modules that are
free as $T$-modules.  If $\Lambda$ has finite representation type, then
$\Lambda_\p$ is nonsingular for all nonmaximal primes $\p$ of $T$
\cite[Theorem 10]{Auslander:isolsing}.  We are primarily interested in
the commutative case, so we leave to the interested reader the extension
of Corollary~\ref{maincor} to Auslander's context.  

We are grateful to the anonymous referee for pointing out this
generalization, as well as a small change in the proof of
Theorem~\ref{main} that makes the extension possible.\end{rem}

An almost immediate consequence of the Theorem~\ref{main} is a rather
powerful
one:

\begin{thm}\label{powerful} Let $(R,\m)$ be a Cohen--Macaulay local ring
and let
$M$ be a finitely generated maximal Cohen--Macaulay $R$-module. Set $e =
e(R)$, the multiplicity of $R$, and $n = \mu(M)$,
the minimal number of generators of $M$. If
$R$ has only finitely many indecomposable nonisomorphic MCM
modules of multiplicity at most $n\cdot e$, then $M_\p$ is free for all
primes
$\p\ne \m$.
\end{thm}

\begin{proof} Consider the exact sequence
$$0\rightarrow N\rightarrow R^n\rightarrow M\rightarrow 0.$$
This is an element of $\Ext^1_R(M,N)$. Since the sum of the
multiplicities of $M$ and
$N$ is $n\cdot e$, Theorem~\ref{main} gives us that $\Ext^1_R(M,N)$ has
finite
length.
Hence after localizing the short exact sequence at an arbitrary prime
$\p\ne \m$ the sequence splits and $M_\p$ is then free.
\hfill\qed\end{proof}

Using this we can give stronger results concerning when a given local
ring is regular, or Gorenstein, or a complete intersection on the
punctured spectrum in terms of the ring having only finitely many
indecomposable MCM modules up to a certain multiplicity.

\begin{cor} Let $(R,\m)$ be either a localization of a finitely
generated algebra over a field $k$ of characteristic $0$, or a quotient
of a power series ring over a field $k$ of characteristic $0$. In the
first case,
let $\Omega$ be the module of K\"ahler differentials of $R$ over $k$,
and in the second
case let $\Omega$ be the universally finite module of differentials of
$R$ over $k$. Assume that $\Omega$ is a maximal Cohen--Macaulay module.
Set $e = e(R)$ and let $n$ be the embedding dimension of $R$, that is,
$n = \mu(\m/\m^2)$. If
$R$ has only finitely many indecomposable nonisomorphic MCM
modules of multiplicity at most $n\cdot e$, then $R_\p$ is regular for
all
primes $\p\ne \m$.
\end{cor}

\begin{proof} In either case there is an exact sequence, 
$$0\rightarrow N\rightarrow R^n\rightarrow \Omega\rightarrow 0.$$
(The number of generators of $\Omega$ is at most the embedding dimension
of
$R$ by Theorem~25.2 of \cite{Matsumura} or Theorem~11.10 of
\cite{Kunz:1986}.)
Applying Theorem~\ref{powerful} then shows that $\Omega_\p$ is a free
$R_\p$
module 
for all primes $\p\ne \m$. Applying Theorem 7.2 of \cite{Kunz:1986} in
the
first case and Theorem 14.1 of \cite{Kunz:1986} in the second case, we
see that $R_\p$ is regular.\hfill\qed
\end{proof}

Although the assumption that the module of differentials is
Cohen--Macaulay
is quite strong, it can occur even if the ring is not regular:  Let $B =
k[X_{ij}]$, where $X =
(X_{ij})$ is an $m\times n$ matrix of
indeterminates, and let $r < {\rm min}(m,n)$ be an integer.  Let $R =
B/I_{r+1}(X)$.  If $I_{r+1}(X)$ has grade at most two, then the module
of
K\"ahler differentials $\Omega$ is a MCM $R$-module
(\cite[Prop. 14.7]{Bruns-Vetter:1988}).

\begin{cor} Let $(R,\m)$ be a Cohen--Macaulay local ring with canonical
module
$\omega$.
Assume that $R$ has only finitely many indecomposable nonisomorphic MCM
modules of multiplicity at most $\type(R)\cdot e(R)$.  Then $R$ is
Gorenstein
on the
punctured
spectrum.\end{cor} 

\begin{proof} Applying Theorem~\ref{powerful} to the
extension 
$$\chi: 0 \to \syz^R_1(\omega) \to R^t \to \omega \to 0,$$
where $t = \type(R)$, shows that $\chi_\p$ is split for every
nonmaximal prime $\p$. It follows that $\omega_{\p}$ is free for all
nonmaximal primes $\p$, which implies that $R_{\p}$ is Gorenstein for
all such primes. \hfill\qed\end{proof}

\begin{cor} Let $R=A/I$, where $A$ is a regular local ring.  Assume
that $R$
is Cohen--Macaulay and that $I/I^2$ is a MCM $R$-module.  If $R$ has
only finitely
many
indecomposable nonisomorphic MCM modules of rank at most $\mu_A(I)$,
then
$R_\p$ is a complete intersection for all nonmaximal primes $\p$ of
$R$.\end{cor} 

\begin{proof} This follows from Theorem~\ref{powerful} as well, together
with
the
fact 
that
$I/I^2$ is a free $R$-module if and only if $I$ is generated by an
$A$-regular sequence (\cite[19.9]{Matsumura}).\hfill\qed\end{proof}

This raises the question of when $I/I^2$ is a MCM $R$-module.  Herzog
\cite{Herzog:1978a} showed that this is the case when $I$ is a
codimension
three prime ideal in a regular local ring $A$ such that $R = A/I$ is
Gorenstein, and this was generalized in \cite{Huneke-Ulrich:1989} to the
case
in which $I$ is licci and generically a complete intersection in a
regular
local ring $A$ such that $A/I$ is
Gorenstein. Since the height of the defining ideal of the non-complete
intersection locus of a Gorenstein licci ideal is known to be bounded,
the above Corollary shows that Gorenstein licci algebras in general
can never have finite CM type, nor even have only finitely many
isomorphism
classes of MCMs of rank at most the embedding codimension of the
algebra. This observation  can be extended to the case in which
$R$ is licci but not Gorenstein by using the module $I/I^2 \otimes
\omega$, which is known to be MCM over $R$ \cite{Bu}.

In a similar vein we can prove:

\begin{cor} Let $R=A/I$, where $A$ is a regular local ring.  Assume
that $R$
is Cohen--Macaulay and normal, and that $I/I^{(2)} $ is a MCM
$R$-module.  If $R$ has
only
finitely many
indecomposable nonisomorphic MCM modules of rank at most $\mu_A(I)$,
then
$R_\p$ is a complete intersection for all nonmaximal primes $\p$ of
$R$.\end{cor}

\begin{proof} As in the proofs above, the assumptions together with
Theorem~\ref{powerful} show that $(I/I^{(2)})_\p $ is $R_\p$-free for
all $\p\ne \m$.
Since $R$ is normal, $R$ satisfies Serre's condition $S_2$, is reduced
and $I_Q$ is generated by a regular sequence for all primes $Q\supseteq
I$
such that the height of $Q/I$ is at most $1$. From the main theorem
of \cite{Smith:2000} it follows that $I_\p$ is generated by a regular
sequence for all $\p\ne \m$. \hfill\qed\end{proof}

We also are able to recover a recent result of Leuschke and Wiegand,
\cite{Leuschke-Wiegand:2000}:

\begin{cor} Let $(R,\m)$ be an excellent local Cohen--Macaulay ring
having finite
CM type. Then the completion $\hat R$ also has finite CM type.
\end{cor}

\begin{proof} In \cite{Wiegand:1998} it was shown that if $R$ is
a Cohen--Macaulay local ring which is of finite CM type and such that
the completion
$\hat
R$ is an isolated singularity, then $\hat R$ has finite CM type.  Since
$R$  
is excellent and has an isolated singularity by Corollary~\ref{maincor},
$\hat R$ also has an isolated 
singularity, and the conclusion follows. \hfill\qed\end{proof}

\section{The $F$-signature}

In this section, let $(R, \m)$ be a Cohen--Macaulay local ring of
dimension
$d$ containing a field of
characteristic $p>0$ and such that $R$ is $F$-finite.  As usual, $q$
will
denote a varying power of
$p$.  We
will consider the direct-sum decomposition 
\begin{equation*}\label{q-decomp}R^{1/q} = R^{a_{1q}} \oplus
M_2^{a_{2q}}
\oplus \cdots \oplus M_{n_q}^{a_{n_q q}}\tag{*}\end{equation*}
of $R^{1/q}$ into indecomposable MCM $R$-modules.  We refer the reader
to \cite{Huneke:CBMS} for basic definitions
concerning the theory of tight closure.  

\begin{defn} The {\em $F$-signature} of $R$ is $s(R) =
\limq\frac{a_{1q}}{q^d}$, provided the limit exists.\end{defn}

Smith and Van den Bergh \cite{Smith-vandenBergh:1997} show that if $R$
is
strongly $F$-regular and the set of indecomposable modules $M_i$
appearing in
the decomposition (\ref{q-decomp}) is finite (as $q$ ranges over all
powers
of $p$), then $\limq\frac{a_{iq}}{q^d}$
exists and is positive for each $i$. For further results in this
direction,
see \cite{Yao}. We will show that $s(R)$ always
exists for Gorenstein local rings $R$, and is positive if and only if
$R$ is
$F$-rational. Recall that a local Noetherian ring of positive
characteristic
$R$ is said to
be \it $F$-rational \rm if ideals generated by
parameters are tightly closed. In the case that $R$ is Gorenstein this
is equivalent to every ideal being tightly closed.  In case $R$ is not
necessarily 
Gorenstein, we can easily prove that if $\limsup\frac{a_{1q}}{q^d} >0$,
then
$R$ is weakly $F$-regular.

For the proofs we use a characterization of tight closure in terms of
Hilbert--Kunz functions. 
Let $(R,\m)$ be a local Noetherian ring of prime
characteristic $p$, and let $I$ be an $\m$-primary ideal. The
\it Hilbert--Kunz function \rm of $I$ is the function taking
an integer $n$ to the length of $R/I^{[p^n]}$, where
$I^{[p^n]}$ is the ideal generated by all the $p^n$th powers of
elements of $I$. The Hilbert--Kunz multiplicity of $I$, denoted
$e_{HK}(I)$ or $e_{HK}(I,R)$, is ${\displaystyle \lim_{q = p^n
\rightarrow \infty}}
\frac{\lambda(R/I^{[q]})}{q^d}$,
where $\lambda(M)$ is the length of $M$. This limit always exists (see,
for
example, Chapter~6 of \cite{Huneke:CBMS}).

We need the following special case of Theorem~8.17 of
\cite{Hochster-Huneke:1990}:  

\begin{thm}\label{hilbert-kunz} Let $(R,\m)$ be a reduced,
equidimensional
complete local ring of prime  
characteristic $p$. Let $I\subseteq J$ be two $\m$-primary ideals. Then
$I^* = J^*$ if and only if $e_{HK}(I) = e_{HK}(J)$. (Here $I^*$ denotes
the
tight closure of $I$.)
\end{thm}

We note that the assumption concerning test elements in \cite[Theorem
8.17]{Hochster-Huneke:1990} 
is automatic in this case since the ring is excellent, reduced, and
local.
See \cite[Theorem 6.1]{Hochster-Huneke:1994}.

\begin{thm}\label{F-ratl_criterion} Assume that $(R, \m)$ is a complete
reduced $F$-finite Cohen--Macaulay
local ring
containing a field of prime characteristic $p$ and let $d = \dim
R$. We adopt the notation from the beginning of this section.
Then 
\newline (1) If $\limsup\frac{a_{1q}}{q^d}>0$, then $R$ is weakly
$F$-regular;
\newline (2) If in addition $R$ is Gorenstein, then $s(R)$ exists, and
is
positive if and only if $R$ is $F$-rational.\end{thm}

\begin{proof} Assume that $R$ is not weakly $F$-regular, that is, not
all
ideals of $R$ are tightly closed. By \cite[Theorem
6.1]{Hochster-Huneke:1994}
$R$ has a test element, and then \cite[Proposition
6.1]{Hochster-Huneke:1990}
shows that the tight closure of an arbitary ideal in $R$ is the
intersection of
$\m$-primary tightly closed ideals. Since $R$ is not weakly $F$-regular,
there
exists an $\m$-primary ideal $I$ with  
$I
\neq I^*$. Choose an element $\Delta$ of $I : \m$ which is not in
$I^*$. Decompose $R^{1/q}$ into indecomposable MCM $R$-modules as in
(\ref{q-decomp}). Then 
\begin{equation*}\begin{split}
\lambda(R/I^{[q]}) - \lambda(R/(I,\Delta)^{[q]}) 
        &=\lambda(R^{1/q}/I R^{1/q}) -
\lambda(R^{1/q}/(I,\Delta)R^{1/q}) \\
        &= a_{1q}\lambda(R/I) + a_{2q} \lambda(M_2/I M_2) + \cdots \\
        & \phantom{blergh} - [ a_{1q} \lambda(R/(I,\Delta)) +
         a_{2q} \lambda(M_2/(I,\Delta)M_2) + \cdots ] \\
        &\geq a_{1q} \lambda(R/I) - a_{1q}  \lambda(R/(I,\Delta)) \\
        &= a_{1q}.
\end{split}\end{equation*}

Dividing by $q^d$ and taking the limit gives on the left-hand side
a difference of Hilbert-Kunz multiplicities,
$$e_{HK}(I,R) - e_{HK}((I,\Delta),R).$$
But by Theorem~\ref{hilbert-kunz}, this difference is zero, showing that
$\limsup\frac{a_{1q}}{q^d}=0$.

Assume now that $R$ is Gorenstein.  To prove (2), it suffices to take $I
= (\xx)$ generated by a system of
parameters and show that the difference $\lambda(M_i/\xx M_i) -
\lambda(M_i/(\xx,\Delta)M_i)$ is zero for all indecomposable nonfree MCM
modules $M_i$. We state this as a separate lemma.

\begin{lem}\label{socle} Let $(R,\m)$ be a Gorenstein local ring and
let
$M$ be an indecomposable nonfree MCM $R$-module. Let
$\xx$ be a  system of parameters for $R$, and let $\Delta \in R$ be a
representative for the socle of $R/(\xx)$.  Then $\Delta M \subseteq \xx
M$.\end{lem} 

\begin{proof} Choose generators $\{m_1, \ldots, m_n\}$ for $M$ and
define a
homomorphism $R \to M^n$ by $1 \mapsto (m_1, \ldots, m_n)$.  Let $N$ be
the cokernel and set $I = \Ann(M)$, so that we have an exact sequence
$$ 0 \to R/I \to M^n \to N \to 0.$$
First assume that $N$ is MCM. Then $\xx$ is regular on $N$ and on $R/I$.
If
$I \subseteq (\xx)$, then $I = I \cap (\xx) = I\xx$, so
$I=(0)$ by NAK. This gives the exact sequence 
$0 \to R \to M^n \to N \to 0.$
Since $R$ is Gorenstein and $N$ is MCM, this sequence must split,
contradicting the fact that $M^n$ has no free summands. So $I
\not\subseteq
(\xx)$. When we kill $\xx$, therefore, the map $\overline R \to
\overline
M^n$ has nonzero kernel, which must contain $\Delta$. Since the elements
$m_1, \ldots, m_n$
generate $M$, this says precisely that $\Delta M \subseteq
\xx M$. 

Finally assume that $N$ is not MCM. Then when we kill $\xx$, there is a
nonzero $\Tor$: 
$$0 \to \Tor_1^R(N,R/(\xx)) \to \overline R \to \overline M^n \to
\overline N \to 0.$$
Again, the map $\overline R \to \overline M^n$ has nonzero kernel, so
$\Delta
\mapsto 0$.
\hfill\qed\end{proof}

Returning to the proof of (2), we have 
$$\lambda(R/(\xx)^{[q]}) - \lambda(R/(\xx,\Delta)^{[q]}) = a_{1q}$$
and 
$$e_{HK}(\xx,R) - e_{HK}((\xx,\Delta),R) = s(R).$$
This shows that the limit exists, and is positive if and only if $\xx$
can be
chosen
to generate a tightly closed ideal, if and only if $R$ is $F$-rational. 
\hfill\qed\end{proof}

We now derive some basic properties of the $F$-signature.

\begin{prop} Let $R$ be a reduced $F$-finite local ring of prime
characteristic $p$, and assume that the $F$-signature $s(R)$
exists. Then
\newline (1) $s(R_\p) \geq s(R)$ for every prime ideal $\p$ of R;
\newline (2) $s(R) = s(\widehat{R})$, where $\widehat{R}$ is the
completion
of $R$.
\end{prop}

\begin{proof} The proposition is clear from the facts that
$(R_\p)^{1/q} \cong
(R^{1/q})_\p$ and $(\widehat{R})^{1/q} \cong \widehat{R^{1/q}}$.
\hfill\qed\end{proof}

\begin{prop}\label{inequality} Let $(R,\m)$ be a reduced $F$-finite
Cohen--Macaulay local ring with infinite residue field, such that the
$F$-signature $s=s(R)$ exists. Then 
$$(e-1)(1-s) \geq e_{HK}(R) -1,$$
where $e = e(R)$ and $e_{HK}(R)$ are the Hilbert--Samuel and
Hilbert--Kunz
multiplicities, respectively.\end{prop}

\begin{proof} Let $\xx$ be a system of parameters generating a minimal
reduction of $\m$, and take a 
composition series
$$(\xx)=I_{e-1} \subset I_{e-2} \subset \cdots \subset I_1=\m$$
with successive quotients isomorphic to $R/\m$.  Then the proof of
Theorem~\ref{F-ratl_criterion} shows that 
$$\lambda(R/I_{j+1}^{[q]})-\lambda(R/I_{j}^{[q]}) \geq a_{1q}$$ 
for each $j = 1, \ldots, e-2$. Dividing both sides by $q^d$, taking the
limit, and adding the inequalities, we obtain
$$e_{HK}((\xx),R) - e_{HK}(\m,R) \geq (e-1)s.$$
Since $(\xx)$ is a minimal reduction for $\m$ and $R$ is Cohen-Macaulay,
 $e_{HK}((\xx),R) =
e((\xx),R)
= e$. This
gives the desired inequality.\qed\end{proof}

\begin{prop}\label{inequality2} Let $(R,\m)$ be a reduced $F$-finite
Cohen--Macaulay local 
ring such that the $F$-signature $s=s(R)$ exists.  Then
$$s(R) \leq \frac{e_{HK}(I)-e_{HK}(J)}{\lambda(J/I)}$$
for every pair of $\m$-primary ideals $I \subseteq J$.  If $R$ is
Gorenstein and has infinite residue field, then equality is attained for
$I$ a minimal reduction of $\m$ and $J=(I,\Delta)$, where $\Delta$
represents a generator of the socle of $R/I$.\end{prop}

\begin{proof} For arbitrary $q=p^e$, decompose $R^{1/q} \cong R^{a_q}
\oplus M_q$, where $M_q$ has no nonzero free direct summands.  Then 
\begin{equation*}\begin{split}
\lambda(R/I^{[q]}) - \lambda(R/J^{[q]}) 
        &=\lambda(R^{1/q}/I R^{1/q}) - \lambda(R^{1/q}/JR^{1/q}) \\
        &= \lambda((R^{a_q}\oplus M_q)/I(R^{a_q}\oplus M_q)) -
\lambda((R^{a_q}\oplus M_q)/J(R^{a_q}\oplus M_q)) \\
        & \geq a_q \lambda(R/I)- a_q \lambda(R/J)\\
        &=  a_{q}\lambda(J/I).
\end{split}\end{equation*}
Rearranging, dividing by $q^{\dim R}$, and taking the limit give the
result.  The statement about equality is just a rewording of the last
sentence of Theorem~\ref{F-ratl_criterion}.\hfill\qed\end{proof}

\begin{cor} Let $(R,\m)$ be a reduced $F$-finite Cohen--Macaulay local
ring
such that the $F$-signature $s=s(R)$ exists. Then $s(R)=1$ if and only
if
$R$ is regular. \end{cor}

\begin{proof} If $R$ is regular, then for all $q = p^e$, $R^{1/q}$ is a
free $R$-module whose
rank is $q^d$. It follows that $s(R) = 1$. Conversely suppose that $s(R)
= 1$.
Then Proposition \ref{inequality} proves that $e_{HK}(R) = 1$. Since
Cohen--Macaulay rings are automatically unmixed, we obtain from
\cite{Wat-Yos} that $R$ is regular. (See also \cite{Hun-Yao} for an
alternate approach.) 
\hfill\qed\end{proof}

\smallskip

We next compute the $F$-signature for several
examples.  It is evident from these examples  that the $F$-signature
gives very delicate information concerning the nature of
the singularity.

\begin{eg}Let $R = k[[x^n, x^{n-1}y, \ldots, y^n]]$, the
$n^{{\text{th}}}$ Veronese
subring of $k[[x,y]]$, where $k$ is a perfect field of positive
characteristic $p$.  Assume that $n\geq 2$ and $p\!\!\!\not| n$.  Then
Herzog \cite{Herzog:1978}
has
shown that $R$ has finite CM type.  Specifically, the indecomposable
nonfree
MCM $R$-modules are the fractional ideals $I_1 = (x,y)$, $I_2 = (x^2,
xy,
y^2)$, \ldots, $I_{n-1} = (x^{n-1}, x^{n-2}y, \ldots, y^{n-1})$.  For
consistency, denote $R$ also by $I_0$.

We have the following decompositions of $R$, $I_1$, \ldots, $I_{n-1}$ as
modules over the ring $R^p$ of $p^{{\text{th}}}$ powers.  All
congruences are
modulo $n$. This decomposition was done by Seibert \cite{Seibert}.
\begin{align*}
R &= \bigoplus_{\underset{0 \leq a,b <p}{a+b\equiv 0}}R^p x^a y^b \oplus
\bigoplus_{\underset{0 \leq a,b <p}{a+b \equiv -1}}I_1^{[p]}
x^a y^b \oplus \cdots \oplus \bigoplus_{\underset{0 \leq a,b <p}{a+b
\equiv 1}}I_{n-1}^{[p]} x^a y^b \\
I_1 &= \bigoplus_{\underset{0 \leq a,b <p}{a+b\equiv 1}}R^p x^a y^b
\oplus \bigoplus_{\underset{0 \leq a,b <p}{a+b \equiv 0}}I_1^{[p]}
x^a y^b \oplus \cdots \oplus \bigoplus_{\underset{0 \leq a,b <p}{a+b
\equiv 2}}I_{n-1}^{[p]} x^a y^b \\
&\vdots \\
I_k &= \bigoplus_{\underset{0 \leq a,b <p}{a+b\equiv k}}R^p x^a y^b
\oplus
\bigoplus_{\underset{0 \leq a,b <p}{a+b \equiv k-1}}I_1^{[p]} 
x^a y^b \oplus \cdots \oplus \bigoplus_{\underset{0 \leq a,b <p}{a+b
\equiv n-1-k}}I_{n-1}^{[p]} x^a y^b \\ 
&\vdots\\
I_{n-1} &= \bigoplus_{\underset{0 \leq a,b <p}{a+b\equiv -1}}R^p x^a y^b
\oplus \bigoplus_{\underset{0 \leq a,b <p}{a+b \equiv -2}}I_1^{[p]}
x^a y^b \oplus \cdots \oplus \bigoplus_{\underset{0 \leq a,b <p}{a+b
\equiv 0}}I_{n-1}^{[p]} x^a y^b
\end{align*}
For any given pair $(p,n)$ and $k < n$, it is easy to compute the number
of
pairs $(a,b)$ with $0 \leq a,b <p$ and $a+b \equiv k$.  For a general
estimate, let $m_k$ be this
number.  Since we have
$|m_k-m_{k-1}| \leq 1$ (even when the indices are taken modulo $n$) and
$\sum_{k}m_k = p^2$, we obtain
$$\left\lfloor \frac{p^2}{n}\right\rfloor \leq m_k \leq \left\lfloor
\frac{p^2}{n}\right\rfloor +r,$$
where $r$ is the remainder upon writing $p^2=Ln+r$, $0 \leq r < n$.  

This implies that if the data from the decompositions is arranged in a
matrix
$A = (a_{ij})$, where $a_{ij}$ is the number of
direct summands
isomorphic to $I_j^{[p]}$ occurring in the decomposition of $I_i$, then
each
entry of $A$ is of the form $(1/n)(p^2 \pm c)$. For any $s \geq 1,$ the
entries of $A^s$ are again $(1/n)$ times a polynomial in $p^{2}$, of
degree $2s$.  This gives $s(R) = 1/n$.

For a concrete example, take $p=5$ and $n=3$.  Then the matrix A is
\begin{equation*} 
\left[\begin{matrix} 8 & 8 & 9 \\ 9 & 8 &8 \\ 8 & 9 & 8
\end{matrix}\right] 
= \frac{1}{3} \left[\begin{matrix} 5^2-1 & 5^2-1 & 5^2 +2 \\ 5^2 +2 &
5^2-1 &
5^2-1 \\
5^2-1 & 5^2+2 & 5^2-1 \end{matrix}\right].
\end{equation*}

The entries of $A^s$ for $s \geq 1$ are polynomials in $p^2=25$ of the
form
$\frac{1}{3}p^{2s} + \text{lower order terms}$.  Thus we see that
$\displaystyle\lim_{s \to
\infty} \frac{a_{1q}}{p^{2s}} = \frac{1}{3}$.\end{eg}

\smallskip

\begin{eg} (See \cite[Thm. 5.4]{Wat-Yos}.)  Let $(R,\m)$ be a
two-dimensional Gorenstein
complete local ring of characteristic $p$. Assume that $R$ is $F$-finite
and $F$-rational. Then $R$
is a double point and is isomorphic to $k[[x,y,z]]/(f)$, where $f$ 
is one of the following:  
\medskip

\begin{tabular}{lllll}
type & equation & char $R$ & $s(R)$ \\
($A_n$) & $f = xy+z^{n+1}$ & $p\geq 2$ & $1/(n+1)$ & ($n \geq 1$)\\
($D_n$) & $f = x^2+yz^2+y^{n-1}$ & $p \geq 3$ & $1/4(n-2)$ & ($n\geq4$)
\\
($E_6$) & $f = x^2+y^3+z^4$ & $p\geq 5$ & $1/24$ \\
($E_7$) & $f = x^2+y^3+yz^3$ & $p \geq 5$ & $1/48$ \\
($E_8$) & $f = x^2+y^3+z^5$ & $p\geq 7$ & $1/120$ 
\end{tabular}

\medskip
We compute the $F$-signature in each of these examples as follows:
In each example a minimal reduction $J$ of the maximal ideal $\m$
has the property that $\m/J$ is a vector space of dimension $1$. Hence
$e_{HK}(J) - e_{HK}(R) = s(R)$ as in the proof of Theorem
\ref{F-ratl_criterion}.
Since $J$ is generated by a regular sequence and is a reduction of
$\m$, $e_{HK}(J) = e(J) = e(\m) = 2$. On the other hand, \cite[Thm.
5.4]{Wat-Yos}
gives the Hilbert-Kunz multiplicity for each of these examples, giving 
our statement (see also \cite{Wat-Yos2}). 
\end{eg}

\medskip

Given that each of the examples above is an invariant ring of a
polynomial ring $S$ under a finite group $G$, and in each case $s(R)=
|G|^{-1}$, one may ask if this is true in general.  The simple example
$R= k[x_1,x_2,x_3]^{S_3}$, where $S_3$ acts naturally by permuting the
variables, shows that something more is needed; $R$ is regular, so
$s(R)=1$.  When $R$ is Gorenstein, we can completely analyze the
$F$-signature in terms of the $R$-module structure of $S$.

\begin{prop}\label{invariants} Let $(R,\m) \subseteq (S,\n)$ be a
module-finite extension of CM local rings of characteristic $p$ with $R$
Gorenstein.  Assume that $R$ is a direct summand of $S$ (as an
$R$-module), that $R/\m = S/\n$ is infinite, and that $s(R)$ and $s(S)$
both exist.  Let $r=\rank_R S$ and let $f$ be the number of $R$-free
direct summands in $S$.  Then $$s(R) \geq f\frac{s(S)}{r}.$$  If in
addition $S$ is regular, then equality holds.\end{prop}

\begin{proof} Let $I$ be a minimal reduction of $\m$, and let
$J=(I,\Delta)$ where $\Delta$ represents the generator of the socle of
$R/I$.  Then by Proposition~\ref{inequality2}, $s(R) = e_{HK}(I,R) -
e_{HK}(J,R)$.  By \cite[2.7]{Wat-Yos}, $e_{HK}(I,R) = e_{HK}(IS,S)/r$
and similarly for $J$.  Proposition~\ref{inequality} then shows that
$s(R) \geq [s(S) \lambda(JS/IS)]/r$.  Finally, write $S\cong R^f\oplus
M$, where $M$ is a MCM $R$-module with no free summands.  Then
Lemma~\ref{socle} and the definitions of $I$ and $J$ show that
$\lambda(JS/IS)=\lambda(JR^f/IR^f)+\lambda(JM/IM)=f$, as desired.

If $S$ is regular, then we may bypass the use of
Proposition~\ref{inequality} by observing that (\cite[1.4]{Wat-Yos})
$e_{HK}(IS,S) = \lambda(S/IS)$ and similarly for $JS$, so $s(R) =
f\frac{s(S)}{r}$.\hfill\qed\end{proof}

Together with Proposition~\ref{inequality}, Proposition~\ref{invariants}
gives the following interesting application to quotient singularities, a
relationship between  order of the group $G$ such that $R$ can be
represented as $S^G$ and the number of free summands of $S$ as an
$R$-module.

\begin{cor} Let $S$ be an $F$-finite regular local ring of
characteristic $p$, with infinite residue field, and let $G$ be a finite
group acting on $S$ with $p \not|\ |G|$.  Set $R=S^G$ and assume that
$R$ is Gorenstein. Write $S=R^f\oplus M$, where $M$ has no nonzero free
direct summands.  Then 
$$|G| = \frac{f}{s(R)} \geq f \frac{e(R)-1}{e(R) -
e_{HK}(R)}.$$\end{cor}

\providecommand{\bysame}{\leavevmode\hbox to3em{\hrulefill}\thinspace}

\end{document}